 \documentclass[final,1p,times]{elsarticle}

 \usepackage{graphics}

\usepackage{amssymb}
 \usepackage{graphicx}
\usepackage[dvipdfmx]{color}

\usepackage{amssymb}

\journal{Physica D}

\begin{document}

\begin{frontmatter}



\title{Pursuit Fractal Analysis of Time-Series Data}

\author{S.Hasegawa\corref{cor1}\fnref{fn2}}
\ead{android.acmos@gmail.com}
\author{H.Anada\corref{cor1}\fnref{fn1}}
\ead{h-anada@tcu.ac.jp}
\author{S.Kanagawa\corref{cor1}\fnref{f2}}
\ead{skanagawa@tcu.ac.jp}
\address[fn2]{
Tokyo City University Graduate Division Graduate School of Engineering Systems Information Engineering
}\address[fn1]{
Tokyo City University Undergraduate Division Faculty of Knowledge Engineering
}

\begin{abstract}
In this study, we present a method to measure changes over time of fractal dimension. We confirmed that our method can calculate the fractal dimension with the same precision as conventional methods, and tracking performance of our method is higher than that of the conventional methods.
\end{abstract}


\end{frontmatter}



\section{Introduction}

In recent years, the fractal dimension of time series data has been measured in many fields. In various fields such as foreign exchange rates \cite{mandel3} and brain waves, some researchers found the relationship between changes in the data and the fractal dimension, and they have studies to apply the fractal dimension to predict changes of data. The concept of 'fractal' was created in 1977 by Mandelbrot\cite{mandel1}. At first the word 'fractal' was used in the field of mathematics. However the word 'fractal' exceeded the field of mathematics and has been used in various fields recently. The fractal dimension of time-series data has the two definitions. The one is the concept to quantify the topological properties of the reconstructed orbit, and the other is a definition to quantify the complexity. In this paper, we study on the quantification of the complexity of the latter.

We must analyze time-series data more prudently than other kinds of data especially. It is impossible to use the box counting method to calculate the fractal dimension in the case of time-series data. However the Burlaga and Klein method, was proposed in 1986, can calculate in this case \cite{ark1}. The method proposed by Higuchi in 1989 improved the Burlaga and Klein method\cite{higuchi1,higuchi2}. 

Lots of researchers measured a fractal dimension of time-series data using the Higuchi method until now. However it is necessary to analyze the changes over time of the time-series data. The Higuchi method can't detect transition of the fractal dimension of time-series data, because he didn't think changes in the complexity of the data. Therefore, the Higuchi method is not suitable for observing changes over time of the fractal dimension.
 However time-dependent change of the fractal dimension is crucial in order to observe time-series data in various fields.

 We propose a pursuit fractal analysis algorithm of time-series data. The primary objective of this study is to deal with changes over time of the fractal dimension of time-series data. The secondary objective of it is that our method and past methods have a same result of the fractal dimension of time-series data without transition. We can calculate a theoretical value of the fractal dimension for fractional Brownian motion (FBM)\cite{fbm}. First of all, we calculate the fractal dimension of FBM by both our method and past methods to conduct numerical verifications as a method of calculating the fractal dimension. Secondly, we make a time series data by combining a different dimensional FBM, and conduct numerical verifications of pursuit performance by using our method and past methods. We confirm that the fractal dimension obtained by our method is close to the theoretical value. 
 Compared with the existing methods, our method is more sensitive to changes of the fractal dimension.

\section{Related Works}
\subsection{Higuchi's method}

Consider the time series data, x(1),x(2),…,x(N), to be analyzed. 
First, they construct a new time series $x^{m}_k$ defined as follows:
\begin{eqnarray}
x^{m}_k ; x(m),x(m+k),x(m+2k),...,x(m+[\frac{N-m}{k}]\cdot k) \ \ \ \  (m=1,2,\cdots,k) \label{1}\nonumber ,
\end{eqnarray}
where $[ ]$ denotes Gauss' notation, both k and m are integer, k indicates the discrete time interval and m indicates the initial time value. They construct $x^{m}_k$ for each k and compute the average length using next equation
\begin{eqnarray}
L_m(k)=\Bigl[ \sum^{[(N-m)/k]}_{i=1}|x(m+i\cdot k)-x(m+(i-1)\cdot k)|\frac{ (N-1)}{[\frac{N-m}{k}\cdot k]} \Bigr] /k \label{2}\nonumber ,
\end{eqnarray}
where N is the total length of the data sequence x and $(N-1)/[(N-m)/k]$ is a normalization factor. An average length is computed for all time series having the same scale k, as the mean of the k lengths $L_m(k)$  for m=1,2,...,k. This procedure is repeated for each k ranging from 1 to $k_{max}$, and they obtain an average length for each k. Finally, if $<L_{m}(k)>\propto k^{-D}$ then the time-series data x is fractal with the dimension D.
\subsection{BK method}
Consider the time series data, x(1),x(2),…,x(N), to be analyzed. 
They defined the length $L_{BK}(k)$ of the curve x(t) as
\begin{eqnarray}
L_{BK}(k)=\frac{\sum^{[N/k]}_{j=1}|\overline{{X}_{j+1}(k)}-\overline{{X}_{j}(k)}|}{k} \label{3}\nonumber ,
\end{eqnarray} 
where $\overline{X}_j(t)$ denotes the mean of the time series data $x_{(j-1)\cdot k+1}$ to $x_{(j-1)\cdot k+k}$. If $L_{BK}(k)\propto k^{-D}$ then the data x is fractal with the dimension D.

\section{Pursuit Fractal Analysis}
Consider the time series data, x(1),x(2),…,x(N), to be analyzed. 
The first step, we construct a new time series in the same way as Higuchi method, $x^m_k$  is defined as follows:
\begin{eqnarray}
x^{m}_k ; x(m),x(m+k),x(m+2k),...,x(m+[\frac{N-m}{k}]\cdot k)\ \ \ \  (m=1,2,\cdots,k)\label{4}\nonumber ,
\end{eqnarray}
where $[ ]$ denotes Gauss' notation, both k and m are integer, k indicates the discrete time interval and m indicates the initial time value. When we define the length of each data set, we multiply the weight exponentially as newer data become large. We construct $x^m_k$ for each curve of the time series data, the average length is computed as

\begin{eqnarray}
L_m(k)=\{ \sum^{[(N-m)/k]}_{i=1}(1-\alpha)^{([(N-m)/k]-i)}|x(m+i\cdot k)-x(m+(i-1)\cdot k)|\frac{ (N-1)}{[\frac{N-m}{k}]\cdot k \cdot A} \} /k \label{5}\nonumber ,
\end{eqnarray}
where $\alpha$ is the weight of importance of the new data, A is the sum of each weight, N is the total length of the data x and $(N-1)/([(N-m)/k]\cdot k)$ is a normalization factor. $\alpha$ depends on the number of elements of the data set to take a range [0,1]. An average length is computed for all the data x having the same scale k, as the mean of the k lengths $L_m(k)$ for m=1,2,...,k. This procedure is repeated for each k ranging from 1 to k to obtain an average length for each k. 

Finally, the average value $<L_m(k)>$, of the sets of $L_m(k)$, is as follows
\begin{eqnarray}
<L_m(k)>=\frac{\sum^{k}_{m=1}L_m(k)}{k}\nonumber .
\end{eqnarray} 
If $<L_m(k)>\propto k^{-D}$ then the time series data x is fractal with the dimension D.

\section{Validation Methodology and Result}
We confirm efficiency of our method from two directions.

First, we calculate the fractal dimension of FBM with the theoretical value 1.5 by three methods, the Higuchi's method, the BK method and our method. The results are shown in Table \ref{table1}.

\begin{table}[htb]
 \centering
 \caption{Result of application to FBM data.}
 \small
  \begin{tabular}{|l||c|c|c|} \hline
    Number of Data & n=2\verb|^|15 & n=2\verb|^|16 & n=2\verb|^|17 \\ \hline
    BK Method & 1.521644$\pm$0.0002896329 & 1.518476$\pm$0.0002185431 & 1.515751$\pm$ 0.0001691493 \\
    Higuchi Method & 1.500634$\pm$0.0001690080&1.500573$\pm$0.0001273236& 1.500378$\pm$ 9.805165e-05 \\
    Our Method  & 1.501013$\pm$0.0002735483&1.501039$\pm$0.0002115390 &1.500559$\pm$0.0001607777 \\ \hline
  \end{tabular}
 \label{table1}
\end{table}
As Table \ref{table1} shows, our method is close to the theoretical value as the Higuchi method. Variance of our method and the BK method are close.

Secondly, we synthesize the FBM with the theoretical fractal dimension 1.5 and the FBM with the theoretical fractal dimension 1.3 as follows:
\begin{eqnarray}
x_{i}(t):x_{1.3}(t_n),x_{1.3}(t_{n-1}),...,x_{1.3}(t_3),x_{1.3}(t_2),0,x_{1.5}(t_2),x_{1.5}(t_3),...,x_{1.5}(t_n)\nonumber ,
\end{eqnarray}
where t indicates the time, and i indicates the fractal dimension.
Then, we examine which method follows the changes of data more quickly.
We calculated the four cases for $\alpha$ based on the number n of components of the dataset($\alpha = 1/n, 1/0.5n, 1/0.3n, 1/0.25n$).

\begin{figure}[t]
\includegraphics[scale=0.6]{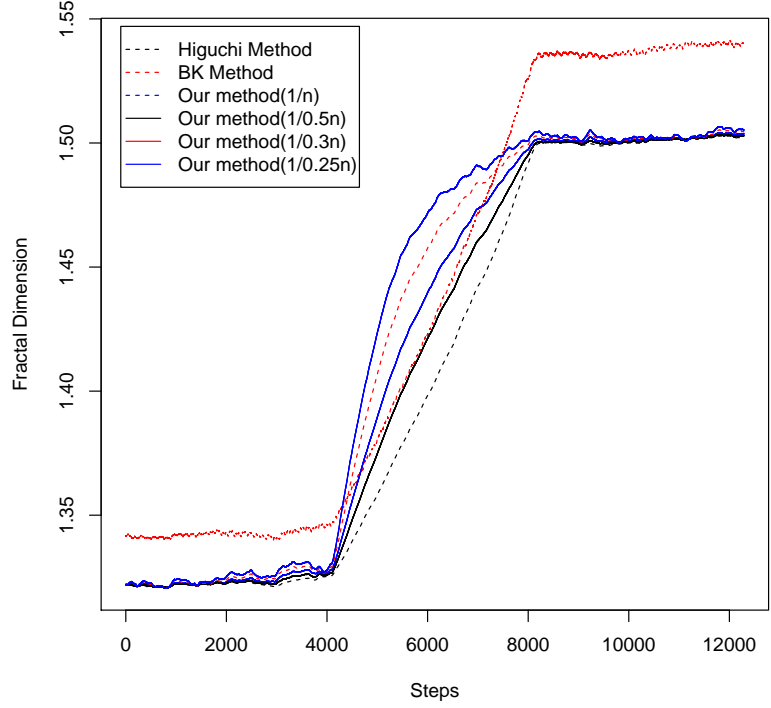}
\centering
\caption{Fractal Dimension :
In this figure, the black broken line shows Higuchi method, the red broken line shows BK method, the blue broken line shows our method when $\alpha =(n)^{-1}$, the black solid line shows our method when $\alpha =(0.5n)^{-1}$, the red solid line shows our method when $\alpha =(0.3n)^{-1}$ and the blue solid line shows our method when $\alpha =(0.25n)^{-1}$. Where n is the number of components of the dataset. The number of trial is 100. We synthesize two sets of time series data by adding the data with the theoretical fractal dimension 1.5 (number of data $2^{13}$) after the data with the theoretical fractal dimension 1.3 (number of data $2^{13}$). Length of analysis range is $2^{12}$ and are shifted one by one from the beginning to the end to analyze.}
\label{fig1} \end{figure}

We calculate fractal dimensions of different subsets, $2^{12}$ data, in the full data, $2^{14}$ data, set from the top of the data to the bottom of the data shifting one by one. The average of the results of 100 times independent repetitions of a calculation experiment is shown in Figure \ref{fig1}.

 The theoretical value is 1.3 until 4096 steps and 1.5 after 4096 steps. It is desirable that our method converge to 1.5 more quickly than the other two methods after 4096 steps.

Figure\ref{fig1} shows that tracking performance of our method is higher than the that of the other two methods. It has a marked influence effect on tracking performance that the weighting for each newer data increases exponentially. 
However, the larger $\alpha$, the more unstable the fractal dimension calculated by our method becomes.
 The value calculated by other methods is stable after about 8000 steps.

\section{Conclusion and Discussion}

We propose a method to measure the changes over time of fractal dimension. We confirmed that our method can calculate the fractal dimension with the same precision as conventional methods, and the tracking performance of our method is higher than that of the conventional methods.

The first result, the fractal dimension obtained by our method is close to the theoretical value as the Higuchi method. Variance of the fractal dimension obtained by our method is close to the one obtained by the BK method.
Our method is considered to have a loss of precision, in the case of calculating the fractal dimension because our method makes exponentially-weighted time-series data. However, the average value is almost the same accuracy as the Higuchi method, the variance is almost the same size as the BK method. We consider that there is no major problem in our method to calculate the fractal dimension of time-series data.
 The second result, there are two discussion points. The first one, the value of the BK method of Figure 1 is away from the theoretical value. Because the data set is too small for the BK method to obtain the fractal dimension. In the case of using less data by the generation method of FBM in this study, the fractal dimension by the BK method tends to be higher values. However, the purpose of this study is to construct a method to measure the changes over time of fractal dimension. It is crucial whether our method can measure the changes over time of fractal dimension enough or not.
The secound, our method has a highest tracking accuracy compared to conventional methods.
We multiply the weight exponentially as newer data become large in our method. 
$\alpha$ is a parameter to adjust the weight of the exponential function. We think to change the $\alpha$ depending on the length of the number of components of the dataset.
When the $\alpha$ is small, the fractal dimension calculated by our method is stable, but tracking performance is low. On the other hand, when the $\alpha$ is large, the fractal dimension calculated by our method is unstable, but the tracking performance is high.
Our method trade off stability against tracking performance. Then our method is sensitive to small variation in time series data but the value is not stable.

There are two problems in this study. 
First, our method has some limitations, such as a limit to the analysis of time series data without an eminent peak. However most observation data include eminent peak and periodic component. 
It remains a challenge for future research to confirm the effectiveness of our method to analyze observational data.
Second, we know that it is necessary to get more data to increase the accuracy of our analysis. However, we don't know the necessary number of time-series data sampling by our method and we do not know the correlation between accuracy and the number of data. A detailed analysis on verification of efficiency of our method in the case of observation data will be presented in a future paper.

\bibliographystyle{elsarticle-harv}
\bibliography{<your-bib-database>}







\end{document}